\documentclass[11pt]{amsart}

\usepackage{amssymb,latexsym,vmargin}

\newtheorem{theorem}{Theorem}
\newtheorem{corollary}{Corollary}

\theoremstyle{definition}
\newtheorem{definition}{Definition}

\def \a {\alpha}

\def \c {\gamma}

\def \< {\langle}
\def \> {\rangle}

\def \supp {{\rm supp}}

\newcommand{\beql}[1]{\begin{equation}\label{#1}}
\newcommand{\eeq}{\end{equation}}
\newcommand{\comment}[1]{}

\newcommand{\Abs}[1]{{\left|{#1}\right|}}

\newcommand{\Norm}[1]{{\left\|{#1}\right\|}}

\newcommand{\Set}[1]{{\left\{{#1}\right\}}}

\newcommand{\RR}{{\mathbb R}}

\newcommand{\ZZ}{{\mathbb Z}}

\newcommand{\TT}{{\mathbb T}}

\newcommand{\inner}[2]{{\langle #1, #2 \rangle}}

\newcommand{\dens}{{\rm dens\,}}

\newcommand{\vol}{{\rm vol\,}}
\newcommand{\ft}[1]{\widehat{#1}}

\newcounter{rem}
\begin{document}

\title[Covering the plane by rotations]{Covering the plane by rotations of a lattice arrangement of disks}

\author[A. Iosevich]{Alex Iosevich}
\address{A.I.: Department of Mathematics, University of Missouri, Columbia, U.S.A}
\email{iosevich@math.missouri.edu}

\author[M. Kolountzakis]{Mihail N. Kolountzakis}
\address{M.K.: Department of Mathematics, University of Crete, Knossos Ave.,
GR-714 09, Iraklio, Greece}
\email{kolount@gmail.com}

\author[M. Matolcsi]{M\'at\'e Matolcsi}
\address{M.M.: Alfr\'ed R\'enyi Institute of Mathematics,
Hungarian Academy of Sciences POB 127 H-1364 Budapest, Hungary.}
\email{matomate@renyi.hu}

\date{November 2006}

\thanks{
M.K.: Supported by the Greek research program ``Pythagoras 2'' (75\% European funds
and 25\% National funds) and by INTAS 03-51-5070 (2004)
({\em Analytical and Combinatorial Methods in Number Theory and Geometry}).
$\bullet$ M.M.: Supported by Hungarian research funds OTKA-T047276, T049301, PF 64061.
}

\begin{abstract}
Suppose we put an $\epsilon$-disk around each lattice point in the plane, and then
we rotate this object around the origin for a set $\Theta$ of angles.
When do we cover the whole plane, except for a neighborhood of the origin?
This is the problem we study in this paper.
It is very easy to see that if $\Theta = [0,2\pi]$ then we do indeed cover.
The problem becomes more interesting if we try to achieve covering with a small closed
set $\Theta$.
\end{abstract}

\maketitle

\section{Introduction}
\label{sec:intro}

In this paper we discuss problems of covering the plane, or all but a bounded part of it,
by rotations of fattened lattices.

Let $\Lambda \subseteq \RR^2$ be a lattice in the plane (a discrete subgroup of $\RR^2$, of dimension 2)
and $\epsilon>0$ be a small number. We define the fattened lattice
$$
E = E(\Lambda,\epsilon) = \Lambda+B_\epsilon(0),
$$
as the $\epsilon$ neighborhood of $\Lambda$ (here $B_\epsilon(0) = \Set{x\in\RR^2: \Abs{x} < \epsilon}$).

Suppose, as we shall do throughout this paper, that $\Theta$ is a set of angles, viewed
as a subset of $S^1$, the unit circle in the plane.
We shall always assume that $\Theta$ is a closed set (see the remark after Definition 1).
If $R_\theta$ denotes the rotation by $\theta$ and
$$
R_\Theta E = \Set{R_\theta x: \theta\in\Theta, x\in E},
$$
the question we are interested in is when $R_\Theta E$ contains the complement of a disk,
when, in other words, $E$ rotated by the angles in $\Theta$ covers everything except
the only obvious obstacle, a neighborhood of the origin.

It is easy to see, and left to the reader, that if we rotate by all possible angles, namely if we take $\Theta = S^1$, then we do indeed achieve covering. 
The question becomes interesting if we try to achieve the same with a {\em small} closed set $\Theta$.

This problem was motivated by earlier results on distances appearing between points of a set of positive upper density. In fact, a question raised by Sz.\ R\'ev\'esz was whether for any set $E$ of positive upper density, the union of finitely many rotates of $E-E$ can cover the complement of a disk. We answer this question in the negative (Theorem \ref{th:no-finite}).   
The first positive result we obtained in this circle of problems (Corollary \ref{cor:smear}),
was deduced easily using a result (Theorem \ref{th:good-measure})
which speaks about which distances
are realizable in sets of positive upper density in Euclidean spaces. 
Theorem \ref{th:good-measure} was obtained in \cite{kolountzakis}
by a careful rewriting of an earlier result of Bourgain
\cite{bourgain} who had improved on Falconer and Marstrand \cite{falconer-marstrand}
and Furstenberg, Katznelson and Weiss \cite{furstenberg-katznelson-weiss}.

\begin{definition}
The set of angles $\Theta \subseteq S^1$ is called $(\Lambda,\epsilon)$-{\em good}
if $R_\Theta E$ contains the complement of a disk, where $E = \Lambda+B_\epsilon(0)$.
The set $\Theta$ will be called {\em good} if it is $(\Lambda,\epsilon)$-good for all
lattices $\Lambda$ and $\epsilon>0$.
\end{definition}

It is easy to see that $\Theta$ is $(\Lambda,\epsilon)$-{\em good} if and only if its closure $\overline\Theta$ is $(\Lambda,\epsilon)$-{\em good}. Therefore, we restrict our attention to closed sets throughout this paper.

In summary our results are as follows.
\begin{enumerate}
\item
If $\Theta\subseteq S^1$ is any arc then $\Theta$ is good (Corollary \ref{cor:smear}).
This follows from using Theorem \ref{th:good-measure} which was proved in \cite{kolountzakis}.
We also give an elementary proof of Corollary \ref{cor:smear} in \S\ref{sec:elementary}.
\item
Using Corollary \ref{cor:smear} we prove in Corollary \ref{cor:very-good-sequence}
that for any $\Lambda$ there are sets $\Theta \subseteq S^1$, which consist
of a convergent sequence of angles plus its limit point, and which are
$(\Lambda,\epsilon)$-good for all positive $\epsilon$.
\item
If $\Theta\subseteq S^1$ is finite then $\Theta$ is not $(\Lambda,\epsilon)$-good
for any lattice $\Lambda$ and any $\epsilon$ smaller than half the shortest
non-zero vector in $\Lambda$ (Theorem \ref{th:no-finite}).
\item
For any lattice $\Lambda$ and any $\epsilon$ which is smaller than half the shortest
non-zero vector in $\Lambda$ there exists an infinite closed set $\Theta\subseteq S^1$
which is not $(\Lambda,\epsilon)$-good
(Corollary \ref{cor:bad-sequence}).
We also prove that this set $\Theta$ may be taken to be a perfect (Cantor-type) set
(Corollary \ref{cor:bad-perfect-set}).
\item
If $\Theta \subseteq S^1$ is rich enough to support a probability measure
whose Fourier Transform is small near infinity (depending on $\Lambda$ and $\epsilon$)
then $\Theta$ is $(\Lambda,\epsilon)$-good (Theorem \ref{th:nice-measure-implies-covering}).
Since any arc of $S^1$ supports probability measures whose Fourier Transform tends to $0$
this is a new proof of Corollary \ref{cor:smear}.
Theorem \ref{th:nice-measure-implies-covering} is
proved directly and not by appealing to any results on distance sets.
\item
If $\Theta \subseteq S^1$ has positive one-dimensional measure
then it is good (Corollary \ref{cor:positive-measure}).
\item
There are sets $\Theta \subseteq S^1$ of $0$ one-dimensional measure
which are good (Corollary \ref{cor:zero-measure}).
\end{enumerate}

\noindent
{\bf Open problem}:
Let $\epsilon>0$ and $E = \Set{(x,y): x\in\ZZ, y\in\RR} + B_\epsilon(0)$.
Is there a finite set of angles $\theta_1,\ldots,\theta_n$ such that
$$
\bigcup_{j=1}^n R_{\theta_j} E
$$
covers the plane?

One might try to prove that this is not the case by showing that in any such finite
set of rotations of $E$ {\em any} line $y=\alpha x$ which is not parallel to any of the strips
cannot be covered.
This amounts to covering the real line by finitely many dilates of the function
$f(x) = \sum_{n\in\ZZ} \chi_{(-\epsilon,\epsilon)}(x-n)$.
This is indeed possible, for any $\epsilon>0$, so this approach to the open problem above
fails.

\noindent
{\bf Acknowledgment}: We are grateful to Prof.\ Yitzhak Katznelson for showing us the proof of 
Theorem \ref{th:measure-0} as well as that of Corollary \ref{cor:good-sequence}.

\section{Continuous moving}
\label{sec:continuous}

The purpose of this section is to show that any arc is good.
A probability measure is called $\delta$-good if its Fourier Transform is $<\delta$ near infinity.
In \cite{kolountzakis} the following theorem is proved (but not stated in this form).
\begin{theorem}\label{th:good-measure}
Suppose that $E \subseteq \RR^d$, $d\ge 2$, has upper density equal to $\epsilon>0$
and that the $0$-symmetric convex body $K$ affords a $(C_d\epsilon)$-good probability
measure $\sigma$ supported on its boundary (the constant $C_d$ depends on the dimension only).
Then, there exists a nonnegative number $t_0$ such that for all $t\ge t_0$ there exist
$x, y \in E$ with
$$
\Norm{x-y}_K = t\ \ \mbox{and}\ \ \frac{x-y}{\Norm{x-y}_K} \in \supp\sigma.
$$
\end{theorem}

\begin{corollary}\label{cor:smear}
Suppose that $\Lambda = A \ZZ^2 \subseteq \RR^2$ is a lattice and $\epsilon>0$.
Write $E = \Lambda + B_\epsilon(0)$.
Then, for any arc $\Theta \subseteq S^1$ we have
$B_{t_0}^c \subseteq R_\Theta E$,
for some $t_0 > 0$.
\end{corollary}

\begin{proof}
Assume $\Theta = [-\theta_0, \theta_0]$.
Let $\Gamma = [a,b]$ be an arc of $S^1$ of length smaller than $\theta_0$ and
take a smooth probability measure $\sigma$ on $S^1$ whose support is $\Gamma$.
Since $\ft\sigma$ tends to $0$ at $\infty$ we can apply Theorem \ref{th:good-measure}
to the set $E' = \Lambda+B_{\epsilon/2}(0)$ and $\sigma$ and we get that there is $t_0$
such that for any $t \ge t_0$ we have (notice that $E = E'-E'$)
$$
t\Gamma \cap E \neq \emptyset.
$$
This implies that
$$
t\Gamma \subseteq R_{[-\theta_0, \theta_0]} E,\ \ \mbox{for $t\ge t_0$}.
$$
Since finitely many rotations of $\Gamma$ will cover $S^1$, it follows by applying our Theorem \ref{th:good-measure}
finitely many times and taking the maximum $t_0$ that there is a finite $t_0'$ such that
any vector of length $\ge t_0'$ is in $R_{[-\theta_0,\theta_0]}E$.
\end{proof}

\section{Elementary proof of Corollary \ref{cor:smear}}
\label{sec:elementary}

We will give the elementary proof for the lattice $\Lambda =\ZZ^2$ for simplicity. The same idea applies to any other lattice, too.

The covering
$$B_{t_0}^c \subseteq R_{[-\theta_0,\theta_0]}E$$
is clearly equivalent to the fact that each 'annulus-arc' $A_{t, \c}=\{(r,\phi ): \ t<r<t+\epsilon , \c \le\phi\le \c+2\theta_0\}$ (given in polar coordinates) contains a lattice point for any $t>t_0$ and any $\c$.

Take finitely many points $n_j=(\cos \a_j ,\sin \a_j)$, $j=1, \dots N$ on the unit circle such that $\tan \a_j$ is irrational, and every open arc of length $\theta_0$ contains at least one $n_j$. Consider the lines $y=-\frac{1}{\tan \a_j} x$ on the torus $\TT=[0,1]\times [0,1]$. Each of these lines form a dense set on the torus, therefore there exist numbers $h_j>0$ such that the line-segments $S_1= \{(x,y): \ y=-\frac{1}{\tan \a_j}x; \ x\in [0, h_j]\}$ and also $S_2=\{(x,y): \ y=-\frac{1}{\tan \a_j}x; \ x\in [-h_j, 0]\}$ are already $\epsilon /4$ dense in $\TT$ (i.e. for every $q\in \TT$ there is a point $s$ of the segment such that $|s-q|< \epsilon/4$; equivalently, the $\epsilon/4$-neighbourhood of $S_1, S_2$ already covers the whole torus). Let $H=\max \{h_j:  \ j=1, \dots N\}$. It follows, by construction, that for each $j$ the $\epsilon/4$-neighbourhood of {\it any line segment} (i.e. not necessarily starting from the origin) of length $H$ and steepness $-\frac{1}{\tan \a_j}$ covers the whole torus.

Take now any $A_{t, \c}$. There is an $\a_j$ such that $\c <\a_j< \c+2\theta_0$. Consider the point $p$ with polar coordinates $p=(t+\epsilon/2, \a_j)\in A_{t, \c}$. It is clear from plane geometry that if $t$ is large enough then there there is a strip $S$ of steepness $-\frac{1}{\tan \a_j}$ and half-width $\epsilon/4$ and length $H$, starting from $p$ (in one of the directions along the line with steepness $-\frac{1}{\tan \a_j}$), which remains fully inside $A_{t, \c}$. By construction, this strip covers the whole torus, and hence contains a lattice point.

\section{Covering using a convergent sequence of rotation angles}
\label{sec:convergent}

The following is a consequence of Corollary \ref{cor:smear} which was shown to us by Y. Katznelson.
\begin{corollary}
\label{cor:good-sequence}
Suppose that $\Lambda\subseteq \RR^2$ is a lattice, $\epsilon>0$.
Let $I$ be any arc in $S^1$.
We can find a convergent sequence of angles $\theta_n \in I$, $n=1,2,\ldots$,
such that the set $\Theta = \overline{\Set{\theta_n, n=1,2,\ldots}}$ is $(\Lambda,\epsilon)$-good.
\end{corollary}
\begin{proof}
Write $E = \Lambda + B_\epsilon(0)$.
Choose any sequence of arcs $I_n \subseteq I$ which converges to a single point $\theta' \in I$.
From Corollary \ref{cor:smear} there is an increasing sequence of numbers $r_n\to \infty$ such that
$$
B_{r_n}^c \subseteq R_{I_n} E.
$$
Let $F_n$ be a finite subset (by compactness such a subset exists) of $I_n$ such that
$$
\overline{B_{r_{n+1}}} \setminus B_{r_n} \subseteq R_{F_n} E.
$$
It follows that the countable set $F = \bigcup_{n=1}^\infty F_n$ is such that
$$
B_{r_1}^c \subseteq R_F E.
$$
Obviously $F$ is a sequence that converges to $\theta'$.
\end{proof}

Corollary \ref{cor:good-sequence} can be strengthened as follows.
\begin{corollary}
\label{cor:very-good-sequence}
For any lattice $\Lambda\subseteq \RR^2$ we can find a convergent sequence of angles $\theta_n$
such that the set $\Theta = \overline{\Set{\theta_n, n=1,2,\ldots}}$ is $(\Lambda,\epsilon)$-good
for all $\epsilon>0$.
\end{corollary}

\begin{proof}
Pick a positive sequence $a_n \to 0$ and, using Corollary \ref{cor:good-sequence},
find a set $\Theta_n \subseteq (0, a_n)$, which consists of a sequence convergent to $0$,
such that $\Theta_n$ is $(\Lambda,1\bigl/n)$-good.
Clearly the set $\bigcup_{n=1}^\infty \Theta_n$ is a sequence which converges to $0$ and is
$(\Lambda,\epsilon)$-good for all positive $\epsilon$.
\end{proof}

\section{Finitely many rotations are not enough, and so are some infinite sets}
\label{sec:finitely-many}

\begin{theorem}\label{th:no-finite}
Let $\Lambda$ be a lattice in $\RR^d$, $d\ge 2$, and $\epsilon>0$ be smaller than
$s(\Lambda)/2$, where $s(\Lambda)$ is the length of the shortest non-zero vector of $\Lambda$.
Write as usual $E = \Lambda+B_\epsilon(0)$.
Then it is impossible to find a finite set of orthogonal matrices
$O_1,\ldots,O_n$ such that $\bigcup_{j=1}^n O_j E$ contains the complement
of a ball.
\end{theorem}
\begin{proof}
Suppose $B_r(0)^c \subseteq \bigcup_{j=1}^n O_j E$.

Let $\epsilon<\epsilon'<s(\Lambda)/2$ and take
$\phi\ge0$ to be a continuous function with $\supp\phi = B_{\epsilon'}(0)$
which is $\ge 1$ on $B_\epsilon(0)$.
Then the functions $f_j(x) = \sum_{\lambda\in O_j\Lambda} \phi(x-\lambda)$ are
periodic continuous functions and writing $f = \sum_{j=1}^n f_j$ we have
\beql{eq-f}
B_r(0)^c \subseteq \bigcup_{j=1}^n O_j E \subseteq \Set{x\in\RR^d: f(x)\ge 1}.
\eeq
It follows that $f$ is almost-periodic (see \cite[p.\ 59]{besicovitch})
hence there are arbitrarily large vectors $T\in\RR^d$
such that $\Norm{f(x)-f(x-T)}_{L^\infty(\RR^d)} \le 1/2$.
But there is an annular neighborhood of $0$ where $f=0$.
By the almost periodicity of $f$ this implies that there are translates
of this neighborhood arbitrarily far where $f\le 1/2$, and this contradicts \eqref{eq-f}.
\end{proof}

Using Theorem \ref{th:no-finite} we can prove the following.
\begin{corollary}
\label{cor:bad-sequence}
Assume the notations of Theorem \ref{th:no-finite} and let $\Lambda$ and $\epsilon$ be fixed, with
$\epsilon < s(\Lambda)/2$.
Then there is an infinite $\Theta \subseteq S^1$
such that the set $R_\Theta E$ is not $(\Lambda,\epsilon)$-good.
\end{corollary}

\begin{proof}
Our set $\Theta$ will be $\overline{\Set{\theta_1, \theta_2, \ldots}}$,
where $\theta_n$ is a convergent sequence.
We shall construct $\Theta$ inductively and along with it we shall construct regions which are not
covered by $R_\Theta \overline{E}$ ($\overline{E}$ is the closure of $E$).

Let $\theta_1$ be arbitrary and assume that we have already chosen the angles $\theta_1,\ldots,\theta_n$ and
that $R_\Set{\theta_1,\ldots,\theta_n} \overline{E}$ does not meet 
the closed disks $G_1,\ldots,G_n$, which are such that
the center of $G_j$ is at distance $j$ from the origin, at least.

Choose $\theta_{n+1}$ distinct from $\theta_1,\ldots,\theta_n$
but so close to, say, $\theta_n$ that the (closed) set
$R_\Set{\theta_1,\ldots,\theta_n,\theta_{n+1}} \overline{E}$ is still disjoint from
the disks $G_1,\ldots,G_n$.

Let now $G_{n+1}$ be a closed disk, disjoint from the disks $G_1,\ldots,G_n$,
whose center is at distance $n+1$ from the origin, at least,
and which is disjoint from $R_{\theta_1,\ldots,\theta_{n+1}}\overline{E}$.
The existence of this disk follows from Theorem \ref{th:no-finite}.

This construction implies the preservation of all the ``holes'' in $R_\Theta \overline{E}$.
\end{proof}

There are even uncountable sets which are not good for covering.
\begin{corollary}
\label{cor:bad-perfect-set}
Assume the notations of Theorem \ref{th:no-finite} and let $\Lambda$ and $\epsilon$ be fixed, with
$\epsilon < s(\Lambda)/2$.
Then there is a perfect set $\Theta \subseteq S^1$ which is not $(\Lambda,\epsilon)$-good.
\end{corollary}

\begin{proof}
The proof is similar to that of Corollary \ref{cor:bad-sequence}.
We shall construct $\Theta$ as an intersection of sets $\Theta_n$ which are finite unions
of closed arcs and each arc of $\Theta_n$ will contain two arcs of $\Theta_{n+1}$.

Along with the set $\Theta_n$ we shall construct a sequence of disjoint closed disks $G_n$,
whose centers are at distance at least $n$ from the origin, and are such that
$R_{\Theta_n} \overline{E}$ is disjoint from $G_1,\ldots,G_n$.
Suppose we have already constructed $\Theta_n$ with the above property.
Pick two points in each of the arcs that make up $\Theta_n$ and call this finite set $F$.
We know that $R_F \overline{E}$ does not cover the complement of a disk, so it must leave uncovered ``holes''
arbitrarily far away from the origin.
Pick a closed disk in such a hole far away and call it $G_{n+1}$.
Let now $\Theta_{n+1} \subseteq \Theta_n$ consist of one tiny closed arc around each point of $F$,
so tiny that the disks $G_1,\ldots,G_n,G_{n+1}$ are still disjoint
from $R_{\Theta_{n+1}} \overline{E}$.
We also make sure that these two tiny intervals in each $\Theta_n$-interval are disjoint.
Clearly then $R_\Theta \overline{E}$ is disjoint from all $G_n$ and is an uncountable perfect set.
\end{proof}

\section{Covering when carrying ``good'' measures}
\label{sec:good-measures-covering}

\begin{theorem}\label{th:nice-measure-implies-covering}
Assume $\Lambda$ is a lattice in the plane and $\epsilon>0$.
Write $E = \Lambda+B_\epsilon(0)$.
Then there is $0<\delta(\Lambda, \epsilon) \sim \dens\Lambda \cdot \epsilon^2$ (as $\epsilon\to0$) such that
if $\Theta \subseteq S^1$ carries a probability measure
$\sigma$ with
$
\limsup_{\xi\to\infty} \Abs{\ft{\sigma}(\xi)} < \delta(\epsilon)
$
then the set $R_\Theta E$ contains the complement of a disk.
\end{theorem}
\begin{proof}
Let $\phi\ge0$ be a $C^\infty$ function supported in $B_{10}(0)$ satisfying $\phi(0) = \ft{\phi}(0) = 1$
and with $\ft{\phi} \ge 0$.
Write $\phi_r(x) = r^{-2} \phi(x/r)$ which also has integral $1$ and is supported in $B_{10 r}(0)$.
For large $q>0$ define
$$
f(x) = f_q(x) = \phi_q\cdot(\phi_\epsilon*\delta_{\Lambda}),
 \ \ \mbox{where $\delta_{\Lambda} = \sum_{\lambda\in\Lambda} \delta_\lambda$.}
$$
It is sufficient to show that if $\Abs{x}$ is sufficiently large then there is $q>0$ such that
$$
\int f(R_\theta x) \,d\sigma(\theta) > 0,
$$
or, equivalently, that
\beql{ft-side}
\int f\left( \Abs{x} R_{x/\Abs{x}} \theta \right)\, d\sigma(\theta)>0.
\eeq
Evaluating \eqref{ft-side} on the Fourier side and applying a change of variable we can rewrite \eqref{ft-side}
as
\beql{to-prove}
\int \ft{f}(\xi) \ft{\sigma}(\Abs{x} R_{x/\Abs{x}} \xi) \,d\xi > 0.
\eeq
From the definition of $f$ and the Poisson summation formula
$$
\ft{\delta_\Lambda} = \dens{\Lambda} \cdot  \delta_{\Lambda^*}
$$ (where $\Lambda^* = A^{-\top}\ZZ^2$ is the dual lattice) we get
$\ft{f} = \dens{\Lambda} \cdot  \ft{\phi_q}*(\ft{\phi_\epsilon} \cdot \delta_{\Lambda^*})$.

Thus the left hand side of \eqref{to-prove}, apart from a factor $\dens\Lambda$, can be written as
\beql{two-parts}
\sum_{\lambda\in\Lambda^*} \ft{\phi}(\epsilon\lambda) \int \ft{\phi}(q(\xi-\lambda))
  \ft{\sigma} (\Abs{x} R_{x/\Abs{x}} \xi) \,d\xi =
\stackrel{I}{\overbrace{\mbox{(term for $\lambda = 0$)}}} +
\stackrel{II}{\overbrace{\sum_{0 \neq \lambda \in \Lambda^*}\cdots}}
\eeq
Since $q^2 \ft{\phi}(q\xi)$ is an approximate identity, with $x$ fixed and $q \to \infty$ we have
$$
I = \int \ft{\phi}(q\xi) \ft{\sigma}(\Abs{x}R_{x/\Abs{x}} \xi)\,d\xi \sim q^{-2}.
$$
This will be the main term in the right hand
side of \eqref{two-parts}.

Write $m(r) = \sup_{\Abs{z}\ge r} \Abs{\ft{\sigma}(z)}$.
Our assumption is that $\limsup_{r\to\infty}m(r) \le \delta(\epsilon)$.
For $II$ we have
\begin{eqnarray*}
II &=& \int \ft{\sigma}(\Abs{x} R_{x/\Abs{x}} \xi)
    \sum_{0 \neq \lambda \in \Lambda^*} \ft{\phi}(\epsilon\lambda) \ft{\phi}(q(\xi - \lambda)) \,d\xi\\
&\le& \int m(\Abs{x}\Abs{\xi}) \sum_{0 \neq \lambda \in \Lambda^*} \ft{\phi}(\epsilon\lambda) \ft{\phi}(q(\xi - \lambda)) \,d\xi
\end{eqnarray*}
Write $G(\xi) = \sum_{0 \neq \lambda \in \Lambda^*} \ft{\phi}(\epsilon\lambda) \ft{\phi}(q(\xi - \lambda))$
and let $r_0$ be the length of the shortest non-zero vector in $\Lambda^*$ and $B = B_{r_0/2}(0)$.
Then
\begin{eqnarray*}
II &\le& \int m(\Abs{x}\Abs{\xi}) G(\xi) \,d\xi \\
 &\le & \int_B G(\xi)\,d\xi + \int_{B^c} m(\Abs{x}\Abs{\xi}) G(\xi)\,d\xi\\
 &=& I_1 + I_2.
\end{eqnarray*}
To estimate $I_1$ we use the fact that $\ft{\phi}\ge0$ and that the balls $\lambda+B$ are disjoint,
$\lambda\in\Lambda^*$:
\begin{eqnarray*}
I_1 &=& \sum_{0\neq \lambda \in \Lambda^*} \ft{\phi}(\epsilon\lambda) \int_B \ft{\phi}(q(\xi-\lambda))\,d\xi\\
 &\le& \int_{B+(\Lambda^*\setminus\Set{0})} \ft{\phi}(q\xi)\,d\xi\\
 &=& q^{-2} \int_{q(B+(\Lambda^*\setminus\Set{0}))} \ft{\phi}(\eta)\,d\eta\\
 &\le& q^{-2} \int_{(qB)^c} \ft{\phi}(\eta)\,d\eta\\
 &\le& o(q^{-2})\ \ \ \mbox{(by the rapid decay of $\ft{\phi}$)}.
\end{eqnarray*}
Finally, for $I_2$ we use our assumption about $m(\cdot)$ and the estimate
\begin{eqnarray*}
\int G(\xi)\,d\xi
 &\le& \sum_{\lambda\in\Lambda^*} \ft{\phi}(\epsilon\lambda) \int\ft{\phi}(q(\xi-\lambda))\,d\xi\\
 &=& q^{-2} \sum_{\lambda\in\Lambda^*} \ft{\phi}(\epsilon\lambda)\\
 &=& C(\epsilon) q^{-2},
\end{eqnarray*}
where $C(\epsilon) \sim \vol\Lambda \cdot \epsilon^{-2}$ as $\epsilon\to 0$.
This shows that $I_2 \le C(\epsilon) q^{-2} m(\Abs{x}r_0/2)$
and,
if $m(r)$ is smaller than $\delta(\epsilon) := 1/C(\epsilon)$ near infinity,
then there is a value $R>0$ such that $\Abs{x} > R$ implies
that \eqref{to-prove} holds for some large $q$, as $I$ will be the dominant term in \eqref{two-parts}.
\end{proof}

\begin{corollary}
\label{cor:positive-measure}
Suppose $\Theta \subseteq S^1$ is a closed set with positive one-dimensional measure.
Then $\Theta$ is good.
\end{corollary}
\begin{proof}
By Theorem \ref{th:nice-measure-implies-covering} it is enough to construct, for any $\delta>0$,
a probability measure $\mu_\delta$
supported on $\Theta$ whose FT is at most $\delta$ in a neighborhood of $\infty$.

For this let $x$ be a Lebesgue point of $\Theta$ and
let the $x$-centered arc $J\subseteq S^1$ be such that
$\Theta$ has density $> 1-\frac{\delta}{10}$ in $J$.
Let $\phi$ be a nonnegative smooth function supported on $J$ such that
the $L^1$ distance of $\phi$ and $\chi_J$
is bounded by $(\delta\Abs{J})/10$ and $\int\phi=\Abs{J}$.

Define the following probability measures:
$$
\mu = \frac{\chi_J}{\Abs{J}},\ 
\nu = \frac{\phi}{\Abs{J}},\ 
\mu_\delta = \frac{\chi_{\Theta\cap J}}{\Abs{\Theta \cap J}}.
$$
By our choice of $J$ and $\phi$ it is clear that
$$
\Norm{\mu-\mu_\delta} < \frac{\delta}{2},\ \Norm{\mu-\nu} < \frac{\delta}{2},
$$
hence we also have $\Norm{\nu-\mu_\delta} < \delta$.
Since $\ft{\nu}(\xi) \to 0$ as $\xi\to\infty$ it follows that $\mu_\delta$ has
FT which is at most $\delta$ in a neighborhood of $\infty$, as required.
\end{proof}

\section{Existence of good sets of rotations of measure $0$}\label{sec:existence}

We owe the following result to Y. Katznelson.
\begin{theorem}\label{th:measure-0}
For any arc in $S^1$ there exists a set $\Theta$ of one-dimensional measure $0$ contained in
that arc which carries a probability measure $\sigma$ whose Fourier Transform
tends to $0$.
\end{theorem}
\begin{proof}
We shall construct $\sigma$ as a weak limit point of a sequence of probability measures $\mu_n$.
We set $\mu_1$ to be arc-length on the given arc, smoothly cut-off by a positive
function and normalized to be a probability measure.
It follows that $\ft{\mu_1}(\xi) \to 0$ as $\Abs{\xi} \to\infty$.

Suppose we have constructed the measure $\mu_n$ and its support is the union of arcs
$I_1^{(n)}, I_2^{(n)}, \ldots, I_{m_n}^{(n)}$.
Assume $\Abs{I_1^{(n)}} \ge \Abs{I_2^{(n)}} \ge \cdots \ge \Abs{I_{m_n}^{(n)}}$.

The next measure $\mu_{n+1}$ will be equal to $\mu_n$ on the arcs $I_2^{(n)}, \ldots, I_{m_n}^{(n)}$.
In $I_1^{(n)}$ the measure $\mu_n$ will be replaced by a measure which will be supported
by a finite union of sub-arcs of $I_1^{(n)}$, all of them shorter than $I_{m_n}^{(n)}$.

Let $R_n\ge \max\Set{n,R_{n-1}}$ be such that $\Abs{\ft{\mu_n}(\xi)} \le 1/n$ for all $\xi$ with $\Abs{\xi}\ge R_n$.
To get $\mu_{n+1}$ from $\mu_n$ in the arc $I_1^{(n)}$ we subdivide $I_1^{(n)}$ into $N\ge 2$
equal arcs and in each of them, say in $[a,b]$,  we shift all the mass of $\mu_n$
into a smooth positive bump in the arc $[a,c]$, where $c-a = \min\Set{(b-a)/2, \Abs{I_{m_n}^{(n)}}/2}$.
Clearly we can choose $N$ so large that
\beql{dyadic}
\Abs{\ft{\mu_{n+1}}(\xi) - \ft{\mu_{n}}(\xi)} \le 2^{-n}/n,\ \ \ (\Abs{\xi}\le R_n).
\eeq
The reason is that if $N$ is large enough the functions
$e_\xi(x) = e^{2\pi i \inner{\xi}{x}}$, $\Abs{\xi}\le R_n$,
are almost constant for $x$ in each arc $[a,b]$.

The new measure $\mu_{n+1}$ is supported on the finitely many arcs
$$
I_1^{(n+1)} = I_2^{(n)},\ldots,I_{{m_n}-1}^{(n+1)} = I_{m_n}^{(n)}
$$
followed by the new arcs $I_{m_n}^{(n+1)},\ldots,I_{m_{n+1}}^{(n+1)}$ that came from $I_1^{(n)}$.
Its Fourier Transform stil tends to $0$ at $\infty$ as $\mu_{n+1}$ is a finite union of smooth bumps.

From the construction it follows that
\beql{first-to-zero}
\mu_n(I_1^{(n)}) \to 0,
\eeq
and that
\beql{length-to-zero}
\Abs{\supp{\mu_n}} \to 0.
\eeq
The reason is that when passing from $\mu_n$ to $\mu_{n+1}$ the effect on the
list of arcs that constitute the support of the measure (remember that the arcs in this list are
decreasing in length) is that the first element of the list is removed
and several members are added to the end of the list.
To see \eqref{length-to-zero} observe that after $m_n$ steps all
the arcs that make up $\supp{\mu_n}$ will have been removed
and $\Abs{\supp{\mu_{n+m_n}}}$  
will be at most $(1/2)\Abs{\supp{\mu_n}}$.
And to see \eqref{first-to-zero}
notice that after the same number of steps we will have that the measure of
the largest arc wil be at most $(1/N) \mu_n(I_1^{(n)})$.

Suppose now that $R_n \le \Abs{\xi} < R_{n+1}$. By the definition of $R_n$ we have
\beql{eq1}
\Abs{\ft{\mu_n}(\xi)} \le 1/n.
\eeq
Since the measures $\mu_n$ and $\mu_{n+1}$ only differ in $I_1{(n)}$ we have
\beql{eq2}
\Abs{\ft{\mu_{n+1}}(\xi) - \ft{\mu_n}(\xi)} \le \mu_n(I_1^{(n)}) = \mu_{n+1}(I_1^{(n)}).
\eeq
Finally, if $k\ge 1$, applying \eqref{dyadic} repeatedly we obtain
\beql{eq3}
\Abs{\ft{\mu_{n+1+k}}(\xi) - \ft{\mu_{n+1}}(\xi)} \le 2/n.
\eeq
Combining \eqref{eq1}, \eqref{eq2} and \eqref{eq3} we obtain
\beql{eq4}
\Abs{\ft{\mu_{k}}(\xi)} \le \epsilon_n := 2/n + \mu_n(I_1^{(n)}),\ \ (k\ge n).
\eeq

Suppose now that $\sigma$ is a weak limit of a subsequence of $\mu_n$.
We have shown that if $R_n \le \xi < R_{n+1}$ then $\Abs{\ft{\sigma}(\xi)} \le \epsilon_n$.
Since $R_n \to \infty$ and $\epsilon_n\to 0$ (this follows from \eqref{first-to-zero})
we have proved that the Fourier Transform of
$\sigma$ tends to $0$ at infinity.
Finally, the support of $\sigma$ is contained in the support of $\mu_n$ for infinitely
many $n$, and hence it has Lebesgue measure $0$, because of \eqref{length-to-zero}.
\end{proof}

The following is now immediate from Theorem \ref{th:nice-measure-implies-covering} combined
with Theorem \ref{th:measure-0}.
\begin{corollary}\label{cor:zero-measure}
In any arc of $S^1$ one can find a good set $\Theta$ of one-dimensional measure 0.
\end{corollary}


\end{document}